\documentclass[11pt]{article}
\usepackage{amsmath,amssymb,latexsym, color}
\usepackage{theorem}
\newtheorem{theorem}{Theorem}
\newtheorem{lemma}{Lemma}
\newtheorem{corollary}{Corollary}
\newtheorem{proposition}{Proposition}
{\theorembodyfont{\rmfamily} }

\newcommand{\esssup}{\operatornamewithlimits{ess\,sup}}

\title{Discrete Laguerre-Sobolev  expansions:\\
A Cohen type inequality}
\author{A. Pe\~{n}a   \thanks{Both authors partially supported by MICINN of Spain under Grant
MTM2009-12740-C03-03, FEDER funds (EU), and the DGA project E-64
(Spain)}\,\,, M. L. Rezola$^{*}$
 \\ \\ {\small Departamento de Matem\'{a}ticas and IUMA.}\\ {Universidad
de Zaragoza (Spain).}}

\date{}

\begin{document}

\maketitle

\begin{abstract}
C. Markett proved a Cohen type inequality for the classical Laguerre expansions  in the appropriate weighted $L^{p}$ spaces. In this paper, we get a Cohen type  inequality for the Fourier expansions in terms of discrete Laguerre--Sobolev orthonormal polynomials with an arbitrary (finite) number of mass points. So,
we extend the  result  due to B. Xh. Fejzullahu and  F. Marcell\'{a}n.

\end{abstract}

2000MSC: 42C05


\bigskip

Corresponding author: Ana Pe\~{n}a

Departamento de Matem\'aticas. Universidad de Zaragoza

50009-Zaragoza (Spain)

e-mail: anap@unizar.es

Phone: 349761328, Fax: 34 976761338

\newpage

\section {Introduction and notations}

Littlewood  conjectured in 1948 that for any trigonometric polynomial $F_N(x)=\sum_{k=1}^N a_k e^{in_k x}$ where $0< n_1 < n_2 < \cdots <n_N$, $\, N\ge 2$, and $\vert a_k\vert \ge 1$ for $1\le k \le N$, there holds the estimate from below
$$\int_0^{2\pi} \vert F_N(x) dx \vert \ge C \,{\rm log}N$$
where $C$ is an absolute constant (see \cite{H-L}).

Cohen's inequality \cite{C} was the first result on the way to the solution of this conjecture. Later, inequalities of this type have been established in various other contexts, e.g., on compact group (see  \cite{G-S-T}).

 In  \cite{M2} Markett proved such inequalities for classical  orthogonal polynomial expansions  in the appropriate weighted $L^{p}$ spaces,  here in terms of the highest coefficient.
The main purpose of this paper is to extend these results  to discrete Laguerre-Sobolev expansions. More precisely, we obtain such inequalities, in the appropriate weighted $L^{p}$ spaces, for Fourier expansions in terms of  orthonormal polynomials with respect to an inner product of the form
\begin{equation}\label{innerproduct}
\langle p,q \rangle_S = \frac{1}{\Gamma(\alpha +1)}\int_0^{\infty} p(x)q(x) \,
x^{\alpha} e^{-x} \,dx + \sum_{j=0}^N M_j p^{(j)}(0) q^{(j)}(0),
\end{equation}
where $\alpha  > -1$ and  $\, M_j \ge 0, \, j = 0, \dots, N.$ Such  inner products are called of discrete  Sobolev type.

Recently  in \cite{B-P}, the authors Fejzullahu and  Marcell\'{a}n obtained Cohen type inequalities for orthonormal expansions with respect to the above inner product in the case $N=1$, i.e. at most two masses in the discrete part. In this particular case, the authors benefit from the fact that there are explicit formulas for the connection coefficients which appear in the representation  of discrete Laguerre-Sobolev type polynomials in terms of three  standard Laguerre polynomials (see \cite{km}). For a general discrete Laguerre-Sobolev  inner product, we only know  that these coefficients  are a nontrivial solution of a system of $N+1$ equations on $N+2$ unknowns (see \cite{k1}). If the system is solved, we get an intricate expression with which it is difficult to work. Our  contribution in this paper is that we can assure   that there exists limit of these connection coefficients and this is enough for our purpose.

\bigskip

 Let $\{L_n^{\alpha}(x) \}_{n\ge0}$ be the sequence of Laguerre polynomials, orthogonal on $[0, \infty)$ with respect to the probability measure $\displaystyle d\mu(x)=\frac{1}{\Gamma(\alpha+1)}x^{\alpha}e^{-x}\,dx$ where $\alpha > -1$ and normalized by $L_n^{\alpha}(0)=\begin{pmatrix}
n+\alpha \\
n \\
\end{pmatrix}$ . We denote the orthonormal Laguerre polynomial of degree $n$ by
$$l_n^{\alpha}(x)=\frac{L_n^{\alpha}(x)}{\Vert L_n^{\alpha} \Vert }$$
where $\Vert L_n^{\alpha} \Vert ^2=\int_0^{\infty}{L_n^{\alpha}(x)}^2\,d\mu(x)$.

Let $\{Q_n^{\alpha} \}_{n\ge0}$ be the sequence of discrete
Laguerre--Sobolev orthogonal polynomials  with respect to the inner
product (\ref{innerproduct}) and such that $Q_n^{\alpha}(x)$ and $L_n^{\alpha}(x)$ have the same leading coefficient. We denote by  $$q_n^{\alpha}(x)={\langle Q_n^{\alpha},Q_n^{\alpha} \rangle }_S^{-1/2}Q_n^{\alpha}(x)$$  the orthonormal discrete  Laguerre-Sobolev polynomials. From now on, for simplicity we write $Q_n(x)=Q_n^{\alpha}(x)$ and $q_n(x)=q_n^{\alpha}(x)$.

Laguerre expansions have been investigated mainly in the following two sets of weighted Lebesgue spaces, namely in the classical spaces (\cite{A-W}, \cite{Muc})

\begin{equation*}
L_p(x^{\alpha}dx)=\left\{\begin{array} {ll}
\{f; \quad \int_0^{\infty} |f(x)e^{-x/2}|^p x^{\alpha}dx< \infty \}, \quad  & \hbox{if $1 \le  p < \infty ;$} \\
\{f; \quad \displaystyle\esssup_{0<x<\infty} \,|f(x)e^{-x/2}|< \infty \}, \quad  & \hbox{if $p = \infty ,$}
\end{array}
\right.
\end{equation*}
for $\alpha > -1$ as well as in the spaces

\begin{equation*}
L_p(x^{\alpha p/2}dx)=\left\{\begin{array} {ll}
\{f; \quad \int_0^{\infty} |f(x)e^{-x/2} x^{\alpha/2}|^p dx< \infty \}, \,\,   & \hbox{if $1 \le  p < \infty ;$} \\
\{f; \quad \displaystyle\esssup_{0<x<\infty} \,|f(x)e^{-x/2}x^{\alpha/2}|< \infty \}, \,\,  & \hbox{if $p = \infty ,$}
\end{array}
\right.
\end{equation*}
for $\alpha > -\frac{2}{p}$ if $1 \le  p < \infty$ and $\alpha \ge 0$ if $p = \infty$.

In order to unify the two results we are going to prove, we introduce an auxiliary parameter $\beta$ which  means either $\alpha$ or $\alpha p/2$.

We consider the class $S_p^{\beta}$, $1 \le p \le \infty$, defined as
the space of measurable functions $f$ defined on $[0, \infty)$, such that there exits $f^{(k)}(0)$ for $k=0, \dots ,N$ and if $1 \le  p < \infty$

\begin{equation*}
\Vert f \Vert _{S_p^{\beta}}^p=\Vert f \Vert _{L_p(x^{\beta}dx)}^p+ \sum_{i=0}^{N} M_j\,|f^{j}(0)|^p< \infty ,
\end{equation*}
\noindent where
\begin{equation*}
\Vert f \Vert _{L_p(x^{\beta}dx)}^p=\int_0^{\infty} |f(x)e^{-x/2}|^p x^{\beta}dx, \quad 1 \le p <\infty,
\end{equation*}
and if  $p = \infty$
\begin{equation*}
\Vert f \Vert _{S_{\infty}^{\beta}}=\text{max} \{ \Vert f \Vert _{L_{\infty}(x^{\beta}dx)}, |f(0)|, \dots, |f^{(N)}(0)| \}< \infty,
\end{equation*}
where
\begin{equation*}
\Vert f \Vert _{L_{\infty}(x^{\beta}dx)}=\left\{ \begin{array}{ll}
\displaystyle\esssup_{0<x<\infty} \,|f(x)e^{-x/2}|, & \hbox{if $\beta=\alpha ;$} \\
\displaystyle\esssup_{0<x<\infty} \,|f(x)e^{-x/2}x^{\alpha/2}|, & \hbox{if $\beta=\alpha p/2 .$}
\end{array}
\right.
\end{equation*}
(If some $M_j=0$ the corresponding derivative does not  appear in the maximum.)

Let $f\in S_p^{\beta}$, $1 \le p \le \infty$, then the Fourier expansion in terms of orthonormal discrete Laguerre-Sobolev polynomials $\{q_n\}_{n \ge 0}$, is
\begin{equation*}
\sum_{k=0}^{\infty}\hat{f}(k)\,q_k(x)
\end{equation*}
 where $\hat{f}(k)=\langle f,q_k \rangle_S$.

In the following, $[S_p^{\beta}]$ denotes the space of all bounded linear operators $T$ from the space $S_p^{\beta}$ into itself, endowed with the usual operator norm,
\begin{equation*}
\Vert T \Vert_{[S_p^{\beta}]}=\sup_{0 \not=f \in S_p^{\beta}}\frac{\Vert Tf \Vert_{S_p^{\beta}}}{\Vert f \Vert_{S_p^{\beta}}}.
\end{equation*}

 Let $1 \le p \le \infty$. For a family of complex numbers $\{c_{k,n} \}_{k=0}^n, \, n \in \mathbb{N} \cup \{0 \}$, with $|c_{n,n}|>0$ we define the operators $T_n^{\alpha,S}:S_p^{\beta}\rightarrow S_p^{\beta}$ by
\begin{equation*}
T_n^{\alpha,S}(f)= \sum_{k=0}^n \, c_{k,n}\hat{f}(k)\,q_k.
\end{equation*}

 Let us denote $q_0=\frac{4\alpha+4}{2\alpha+1}$ for $\beta=\alpha$ and $q_0=4$ for $\beta=p \alpha/2$, and let $p_0$  be the conjugate of $q_0$, i.e. $1/p_0+1/q_0=1$. Now, we can state our main theorem, which extends the ones given in   \cite{M2} and \cite{B-P}.

\begin{theorem}\label{cohen}
Let $1 \le p \le \infty$. There exists a positive constant $C$, independent of $n$, such that:

\noindent For $\alpha>-1/2$
\begin{equation*}
\Vert T_n^{\alpha,S} \Vert_{[S_p^{\alpha}]}\ge C\,|c_{n,n}|\left\{
\begin{array}{ll}
 n^{\frac{2\alpha+2}{p}-\frac{2\alpha+3}{2}}, & \hbox{if $1 \le p <p_0$;} \\
  (\log (n+1))^\frac{2\alpha+1}{4\alpha+4}, & \hbox{if $p=p_0,\, p=q_0$;} \\
   n^{\frac{2\alpha+1}{2}-\frac{2\alpha+2}{p}}, & \hbox{if $q_0<p \le \infty$.}
    \end{array}
     \right.
\end{equation*}
For $\alpha>-2/p$ if $1 \le p < \infty$ and $\alpha \ge 0$ if $p=\infty$
\begin{equation*}
\Vert T_n^{\alpha,S} \Vert_{[S_p^{p\alpha/2}]}\ge C\,|c_{n,n}|\left\{
 \begin{array}{ll}
  n^{\frac{2}{p}-\frac{3}{2}}, & \hbox{if $1 \le p <p_0$;} \\
   (\log (n+1))^\frac{1}{4}, & \hbox{if $p=p_0,\, p=q_0$;} \\
    n^{\frac{1}{2}-\frac{2}{p}}, & \hbox{if $q_0<p \le \infty$.}
    \end{array}
    \right.
\end{equation*}
\end{theorem}

This theorem  will be proved in   Section 3. In Section 2, we obtain some new results for discrete Laguerre-Sobolev polynomials, which we will use to establish Theorem \ref{cohen}. More concretely, we prove a technical  lemma that  will be used to deduce a Mehler-Heine type formula for Laguerre-Sobolev polynomials and a sharp estimation for their norm in the appropriate weighted $L_p$ spaces.

In the sequel we  use the  following notation, $a_n \sim b_n$  means that there exist positive constants $c_1$ and $c_2$, such that $c_1a_n \le b_n \le c_2a_n$ for $n$ large enough, while $a_n \cong b_n$  means that the sequence $\frac{a_n}{b_n}$ converges to $1$. Throughout the paper, the values of the constants may change from line to line.

\section {Estimates for  discrete Laguerre-Sobolev polynomials}

Consider  the standard Laguerre polynomials $L_n^{\alpha}$  and the Laguerre-Sobolev polynomials $Q_n$ with the same leading coefficient.

Let us recall some properties of Laguerre polynomials for $\alpha > -1$ (see \cite{sz}). The evaluation at $x=0$ of the polynomials $L_n^{\alpha}$ and its successive  derivatives are given by
\begin{equation*}
(L_n^{\alpha})^{(k)}(0)=\frac{(-1)^k \Gamma(n+\alpha+1)}{(n-k)!\,\Gamma(\alpha+k+1)}, \,k \in \mathbb{N}\cup \{0 \},
\end{equation*}
and their $L_2$-norm is
\begin{equation}\label{norma Laguerre}
\Vert L_n^{\alpha} \Vert^2=\frac{1}{\Gamma(\alpha +1)}\int_0^{\infty} (L_n^{\alpha}(x))^2 \,
x^{\alpha} e^{-x} \,dx=\frac{\Gamma(n+\alpha+1)}{n!\,\Gamma(\alpha+1)}.
\end{equation}
As usual, we denote the derivatives of the  $n$th kernels of Laguerre polynomials by
$$K_n^{(k,h)}(x,y)=\displaystyle\frac{\partial^{k+h}}{\partial{x^k}
\partial{y^h}}K_n(x,y)
=\displaystyle \sum _{i=0}^n
\frac{(L_i^{\alpha})^{(k)}(x)(L_i^{\alpha})^{(h)}(y)}{\Vert L_i^{\alpha} \Vert ^2}$$
with $k,h  \in \mathbb{N}\cup \{0\}$ and the convention $K_n^{(0,0)}(x,y)=K_n(x,y)$.

In the next lemma, we obtain an asymptotic estimate for
$Q_n^{(k)}(0)$, that will play  an important role
along this paper.
\begin{lemma}\label{cocientederivadas}
Let $ Q_n $ be the polynomials orthogonal with respect to the inner product (\ref{innerproduct}). Then the following statements hold:
\begin{itemize}
\item[(a)]

\begin{equation*}
 \frac{Q_{n}^{(k)}(0)}{(L_n^{\alpha})^{(k)}(0)}\cong \left\{
                              \begin{array}{ll}
                                \displaystyle{\frac{C_k}{n^{\alpha+2k+1}}}, & \hbox{for k such that $ M_k>0$;} \\
                                \quad C_k, & \hbox{otherwise,}
                              \end{array}
                            \right.
\end{equation*}
where $C_k$ is a nonzero constant independent of $n$.

\item[(b)]
$$\langle Q_n, Q_n \rangle_S \cong \Vert L_n^{\alpha} \Vert ^2 \,.$$
\end{itemize}
\end{lemma}
\textbf{Proof.} If all the masses in the inner product (\ref{innerproduct}) are zero the result is trivial because $Q_n=L_n^{\alpha}$.
We will prove the result by induction concerning the number of positive masses in the inner product (\ref{innerproduct}).

We take the first mass which is positive, namely $M_{j_1}$ ($j_1 \ge 0$),  and consider the sequence of  polynomials $\{Q_{n,1}\}_{n \ge0}$ orthogonal with respect to the inner product

\begin{equation*}
(p,q )_{1} = \frac{1}{\Gamma(\alpha +1)}\int_0^{\infty} p(x)q(x) \,
x^{\alpha} e^{-x} \,dx +  M_{j_1} p^{(j_1)}(0) q^{(j_1)}(0).
\end{equation*}

The Fourier expansion of the polynomial $Q_{n,1}$ in the orthogonal basis $\{L_n^{\alpha}\}_{n\ge0}$ leads to
\begin{equation*}
Q_{n,1}(x) = L_n^{\alpha}(x) - M_{j_1} Q_{n,1}^{(j_1)}(0) K_{n-1}^{(0,j_1)}(x,0) \,.
\end{equation*}
Therefore
\begin{equation}\label{Qn,1(x)}
Q_{n,1}(x) = L_n^{\alpha}(x) - \frac{M_{j_1} (L_n^{\alpha})^{(j_1)}(0)}{1+ M_{j_1} K_{n-1}^{(j_1,j_1)}(0,0)} K_{n-1}^{(0,j_1)}(x,0) \,,
\end{equation}
and
\begin{equation}\label{equivnormas1}
(Q_{n,1},Q_{n,1})_{1} =\Vert L_n^{\alpha} \Vert ^2+ M_{j_1}\, \frac{ (({L_n^{\alpha})}^{(j_1)}(0))^2}{1+ M_{j_1} K_{n-1}^{(j_1,j_1)}(0,0)}  \,.
\end{equation}
These relationships are very well known in the literature of discrete Sobolev type orthogonal polynomials.

 Taking derivatives $k$ times in (\ref{Qn,1(x)})  and evaluating at $x=0$, we obtain

\begin{equation}\label{Qn,1/Ln}
\frac{Q_{n,1}^{(k)}(0)}{(L_n^{\alpha})^{(k)}(0)} =1- \frac{M_{j_1}K_{n-1}^{(k,j_1)}(0,0)}{1+M_{j_1}K_{n-1}^{(j_1,j_1)}(0,0)}\frac{(L_n^{\alpha})^{(j_1)}(0)}{(L_n^{\alpha})^{(k)}(0)}  \,.
\end{equation}
 Applying the Stolz criterion (see, e.g. \cite{k}), we have
\begin{equation}\label{Kn(k,j1)}
\lim_n \frac{K_{n-1}^{(k,j_1)}(0,0)}{n^{\alpha+k+j_1+1}}=\lim_n \frac{(L_{n-1}^{\alpha})^{(k)}(0)(L_{n-1}^{\alpha})^{(j_1)}(0)}{\Vert L_{n-1}^{\alpha} \Vert ^2 (\alpha +k+j_1+1) n^{\alpha+k+j_1}}\not=0,
\end{equation}
 and therefore
\begin{equation}\label{cociente nucleos} \frac{K_{n-1}^{(k,j_1)}(0,0)}{K_{n-1}^{(j_1,j_1)}(0,0)}\frac{(L_n^{\alpha})^{(j_1)}(0)}{(L_n^{\alpha})^{(k)}(0)}\cong \frac{(\alpha+2j_1+1)}{(\alpha+k+j_1+1)}\frac{(L_{n-1}^{\alpha})^{(k)}(0)}{(L_{n-1}^{\alpha})^{(j_1)}(0)}
\frac{(L_n^{\alpha})^{(j_1)}(0)}{(L_n^{\alpha})^{(k)}(0)}$$ $$\cong \frac{\alpha+2j_1+1}{\alpha+k+j_1+1}.
\end{equation}
Thus, from (\ref{Qn,1/Ln}), (\ref{Kn(k,j1)}) and (\ref {cociente nucleos}), we have
\begin{equation*}
\frac{Q_{n,1}^{(j_1)}(0)}{(L_n^{\alpha})^{(j_1)}(0)}=\frac{1}{1+M_{j_1}K_{n-1}^{(j_1,j_1)}(0,0)}\cong \frac{C_{j_1}}{n^{\alpha+2j_1+1}}
\end{equation*}
and  for $k\not=j_1$
\begin{equation*}
\frac{Q_{n,1}^{(k)}(0)}{(L_n^{\alpha})^{(k)}(0)} \cong 1-\frac{\alpha +2j_1+1}{\alpha+k+j_1+1}\not=0.
\end{equation*}
 So, we achieve (a) for $Q_{n,1}$. Besides, taking limits in (\ref{equivnormas1})  and using again the size of derivatives of Laguerre polynomials, we get (b) for the polynomials $Q_{n,1}$.

 If there are no more positive masses, since $Q_{n,1}=Q_n$  we have concluded the proof. Otherwise, suppose that the results (a) and (b) hold for the sequence of polynomials $\{Q_{n,s-1}\}_{n \ge0}$ orthogonal with respect to the inner product
\begin{align*}
(p,q )_{s-1} &= \frac{1}{\Gamma(\alpha +1)}\int_0^{\infty} p(x)q(x) \,
x^{\alpha} e^{-x} \,dx \\ \nonumber
&+  M_{j_1} p^{(j_1)}(0) q^{(j_1)}(0)+ \dots +M_{j_{s-1}} p^{(j_{s-1})}(0) q^{(j_{s-1})}(0),
\end{align*}
where $j_1 <j_2 < \dots <j_{s-1}$ and all these masses are positive.
Now, we have to prove the result for the polynomials $Q_{n,s}$ orthogonal with respect to
\begin{align*}
(p,q )_{s} &= \frac{1}{\Gamma(\alpha +1)}\int_0^{\infty} p(x)q(x) \,
x^{\alpha} e^{-x} \,dx \\ \nonumber
&+  M_{j_1} p^{(j_1)}(0) q^{(j_1)}(0)+ \dots +M_{j_{s}} p^{(j_{s})}(0) q^{(j_{s})}(0),
\end{align*}
where $M_{j_{s}}>0$.
Since $(p,q )_{s}=(p,q )_{s-1}+M_{j_s} p^{(j_s)}(0) q^{(j_s)}(0)$ we can work as before. Then the Fourier expansion of the polynomial $Q_{n,s}$ in the orthogonal basis $\{Q_{n,s-1}\}_{n\ge0}$ leads to
\begin{equation*}
Q_{n,s}(x) = Q_{n,s-1}(x) - M_{j_s} Q_{n,s}^{(j_s)}(0) K_{n-1, s-1}^{(0,j_s)}(x,0) \,,
\end{equation*}
where  $K_{n,s-1}$ denotes the corresponding $n$th kernel for the sequence $\{Q_{n,s-1}\}$ and $$K_{n,s-1}^{(k,h)}(x,y)=
\displaystyle \sum _{i=0}^n
\frac{Q_{i,s-1}^{(k)}(x)Q_{i,s-1}^{(h)}(y)}{( Q_{i,s-1},Q_{i,s-1})_{s-1} }, \quad k,h  \in \mathbb{N}\cup \{0\}.$$
Therefore, in the same way as in (\ref{Qn,1(x)}) and (\ref{equivnormas1}), we get
\begin{equation}\label{Qn,s(x)}
Q_{n,s}(x) = Q_{n,s-1}(x) - \frac{M_{j_s} Q_{n,s-1}^{(j_s)}(0)}{1+ M_{j_s} K_{n-1,s-1}^{(j_s,j_s)}(0,0)} K_{n-1,s-1}^{(0,j_s)}(x,0) \,,
\end{equation}
and
\begin{equation}\label{equivnormas s}
(Q_{n,s},Q_{n,s})_{s} =(Q_{n,s-1},Q_{n,s-1})_{s-1} + M_{j_s}  \frac{ ({Q_{n,s-1}}^{(j_s)}(0))^2}{1+ M_{j_s} K_{n-1,s-1}^{(j_s,j_s)}(0,0)}  \,.
\end{equation}
Taking derivatives $k$ times in (\ref{Qn,s(x)})  and evaluating at $x=0$, we obtain
\begin{equation}\label{Qn,s/Ln}
\frac{Q_{n,s}^{(k)}(0)}{(L_n^{\alpha})^{(k)}(0)} =\frac{Q_{n,s-1}^{(k)}(0)}{(L_n^{\alpha})^{(k)}(0)}\left[1 - \frac{M_{j_s}K_{n-1,s-1}^{(k,j_s)}(0,0)}{1+M_{j_s}K_{n-1,s-1}^{(j_s,j_s)}(0,0)}
\frac{Q_{n,s-1}^{(j_s)}(0)}{Q_{n,s-1}^{(k)}(0)} \right] \,.
\end{equation}
Applying the Stolz criterion and
 the hypotheses (a) and (b) for $\{ Q_{n,s-1} \}_{n\ge0}$, we can obtain
\begin{equation}\label{K_{n-1,s-1}}
 K_{n-1,s-1}^{(k,j_s)}(0,0) \cong \left\{
                              \begin{array}{ll}
                                C_k \, n^{\alpha+k+j_s+1}, & \hbox{if $k \not=j_1,\dots, j_{s-1}$;} \\
                                C_k \, n^{j_s-k}, & \hbox{if $ k=j_1,\dots, j_{s-1}$,}
                              \end{array}
                            \right.
\end{equation}
where $C_k$ is a nonzero constant.
 Indeed, for $k \not=j_1,\dots, j_{s-1}$,

\begin{align}\label{K_{n-1,s-1}^{(k,js}(0,0)}
\nonumber
&  \lim_n \frac{K_{n-1,s-1}^{(k,j_s)}(0,0)}{n^{\alpha+k+j_s+1}}
 =\lim_n \frac{Q_{n-1,s-1}^{(k)}(0)Q_{n-1,s-1}^{(j_s)}(0)}{(Q_{n-1,s-1},Q_{n-1,s-1})_{s-1}(\alpha+k+j_s+1) \,\, n^{\alpha+k+j_s}} \\ \nonumber
&=\lim_n \frac{Q_{n-1,s-1}^{(k)}(0)}{(L_{n-1}^{\alpha})^{(k)}(0)}\lim_n\frac{Q_{n-1,s-1}^{(j_s)}(0)}{(L_{n-1}^{\alpha})^{(j_s)}(0)}\lim_n  \frac{(L_{n-1}^{\alpha})^{(k)}(0)(L_{n-1}^{\alpha})^{(j_s)}(0)}{\Vert L_{n-1}^{\alpha}  \Vert^2 \,(\alpha+k+j_s+1)\, n^{\alpha+k+j_s}},  \\
\end{align}
and, for $k=j_1,\dots, j_{s-1}$,
\begin{align}\label{K_{n-1,s-1}^{(k,js)}(0,0)}
& \lim_n \frac{K_{n-1,s-1}^{(k,j_s)}(0,0)}{n^{j_s-k}}
=\lim_n  \frac{Q_{n-1,s-1}^{(k)}(0)Q_{n-1,s-1}^{(j_s)}(0)}{(Q_{n-1,s-1},Q_{n-1,s-1})_{s-1}(j_s-k) \,\, n^{j_s-k-1}} \\ \nonumber
&=\lim_n \frac{(L_{n-1}^{\alpha})^{(k)}(0)(L_{n-1}^{\alpha})^{(j_s)}(0)}{\Vert L_{n-1}^{\alpha}  \Vert^2(j_s-k) \,\, n^{\alpha+k+j_s}} \lim_n  n^{\alpha+2k+1} \frac{Q_{n-1,s-1}^{(k)}(0)}{(L_{n-1}^{\alpha})^{(k)}(0)} \lim_n \frac{Q_{n-1,s-1}^{(j_s)}(0)}{(L_{n-1}^{\alpha})^{(j_s)}(0)}.
\end{align}
Then, from (\ref{Qn,s/Ln}), (\ref{K_{n-1,s-1}}) and the hypothesis for $Q_{n,s-1}$, we have
\begin{equation*}
\frac{Q_{n,s}^{(j_s)}(0)}{(L_n^{\alpha})^{(j_s)}(0)}=\frac{Q_{n,s-1}^{(j_s)}(0)}{(L_n^{\alpha})^{(j_s)}(0)}
\frac{1}{1+M_{j_s}K_{n-1,s-1}^{(j_s,j_s)}(0,0)}\cong \frac{C_{j_s}}{n^{\alpha+2j_s+1}},
\end{equation*}
with $C_{j_s}$ a nonzero constant. Moreover, for $k\not=j_s$, taking into account (\ref{K_{n-1,s-1}^{(k,js}(0,0)}), (\ref{K_{n-1,s-1}^{(k,js)}(0,0)}) and the hypothesis for $Q_{n,s-1}$, we can deduce
\begin{align}
\nonumber
& \frac{K_{n-1,s-1}^{(k,j_s)}(0,0)}{K_{n-1,s-1}^{(j_s,j_s)}(0,0)}\frac{Q_{n,s-1}^{(j_s)}(0)}{Q_{n,s-1}^{(k)}(0)}=
\frac{K_{n-1,s-1}^{(k,j_s)}(0,0)}{K_{n-1,s-1}^{(j_s,j_s)}(0,0)}\frac{Q_{n,s-1}^{(j_s)}(0)}{(L_{n}^{\alpha})^{(j_s)}(0)}
\frac{(L_{n}^{\alpha})^{(k)}(0)}{Q_{n,s-1}^{(k)}(0)}\frac{(L_{n}^{\alpha})^{(j_s)}(0)}{(L_{n}^{\alpha})^{(k)}(0)}\\ \nonumber
& \cong \frac{(L_{n}^{\alpha})^{(j_s)}(0)}{(L_{n-1}^{\alpha})^{(j_s)}(0)} \frac{(L_{n-1}^{\alpha})^{(k)}(0)}{(L_{n}^{\alpha})^{(k)}(0)}\left\{
\begin{array}{ll}
\frac{\alpha+2j_s+1}{\alpha+k+j_s+1}, & \hbox{if $k\not=j_1,\dots, j_{s-1}$;} \\
\frac{\alpha+2j_s+1}{j_s-k}, & \hbox{if $k=j_1,\dots, j_{s-1}$,}
\end{array}
\right.\\
\nonumber
& \cong \left\{
\begin{array}{ll}
\frac{\alpha+2j_s+1}{\alpha+k+j_s+1}, & \hbox{if $k\not=j_1,\dots, j_{s-1}$;} \\
\frac{\alpha+2j_s+1}{j_s-k}, & \hbox{if $k=j_1,\dots, j_{s-1}$.}
\end{array}
\right.
\end{align}
Thus, taking limits in (\ref{Qn,s/Ln}) and (\ref{equivnormas s}), we get $(a)$  and $ (b)$ for the polynomials $Q_{n,s}$, i.e.

\begin{equation*}
 \frac{Q_{n,s}^{(k)}(0)}{(L_n^{\alpha})^{(k)}(0)}\cong \left\{
                              \begin{array}{ll}
                                \displaystyle{\frac{C_k}{n^{\alpha+2k+1}}}, & \hbox{if $ k=j_1,\dots, j_{s}$;} \\
                                \quad C_k, & \hbox{otherwise,}
                              \end{array}
                            \right.
\end{equation*}
and $$(Q_{n,s},Q_{n,s})_{s} \cong \Vert L_n^{\alpha} \Vert ^2 .$$
Hence the result follows.$\quad \Box$

Observe that the part $(a)$ of  Lemma \ref{cocientederivadas} is also true for the ratio of the corresponding orthonormal polynomials, and therefore there exists
\begin{equation}\label{limitecocientederivadas}
 \lim_n\frac{q_n^{(k)}(0)}{(l_n^{\alpha})^{(k)}(0)}=\left\{
\begin{array}{ll}
0, & \hbox{for $k$ such that $M_k>0$;} \\
C_k\not=0, & \hbox{otherwise.}
 \end{array}
 \right.
\end{equation}

\bigskip

Consider the following representation of the orthonormal polynomials $q_n$ in terms of the orthonormal Laguerre polynomials $l_n^{\alpha}$  (see \cite[Section 9]{k1} )

\begin{equation}\label{formulaKoekoek}
q_n(x)=\sum_{j=0}^{N+1} b_j(n)x^jl_{n-j}^{\alpha+2j}(x).
\end{equation}

For the inner product (\ref{innerproduct}) with $N=1$, the coefficients $b_j(n)$ was explicitly obtained in \cite{km}, and their estimation was essential to obtain the result in \cite{B-P}.

Now in the general case,  using Lemma \ref{cocientederivadas}, we can prove that  there is always limit of the connection coefficients  $b_j(n)$ for an arbitrary $N$.

\begin{lemma}\label{limitecoeficientesKoekoek}
Let $\{b_j(n)\}_{0}^{N+1}$ be the coefficients in formula (\ref{formulaKoekoek}). Then, there exists
 $$ \lim_n b_j(n)=b_j \in \mathbb{R}, \quad j \in \{0, \dots, N+1\}.$$
\noindent Moreover, the first index $j$ such that $b_j\not=0$ corresponds with the first $j$ such that $M_j=0$ in the inner product (\ref{innerproduct}). (We understand  that if all the masses are positive, then  the unique coefficient $b_j$ different from zero is the last one).
\end{lemma}
\textbf{Proof.} Taking derivatives $k$ times in (\ref{formulaKoekoek})  and evaluating at $x=0$, we deduce
\begin{equation}\label{cocientederivadasortonormales}
\frac{q_n^{(k)}(0)}{(l_n^{\alpha})^{(k)}(0)}=\sum_{j=0}^{k} b_j(n) \left(
\begin{array}{c}
 k \\
 j \\
  \end{array}
   \right)
j!\,A_j(k,n), \quad k\in \{0, \dots, N+1\},
\end{equation}
where $A_0(k, n)=1$ and
\begin{equation}\label{expresionA}
A_j(k, n)=\displaystyle \frac{(l_{n-j}^{\alpha+2j})^{(k-j)}(0)}{(l_n^{\alpha})^{(k)}(0)}\cong \frac{(-1)^j\Gamma(\alpha+k+1)}{\Gamma(\alpha+k+j+1)}\left(\frac{\Gamma(\alpha+2j+1)}{\Gamma(\alpha+1)}\right)^{1/2}
\end{equation}

Since there exists $\lim_n A_j(k, n)\not=0$ , applying recursively (\ref{limitecocientederivadas}) and (\ref{cocientederivadasortonormales})  we can assure  there exists $\lim_n b_j(n)=b_j, \, j \in \{0, \dots, N+1 \}$.
More precisely,  for $k=0$ we have
$$\lim_n b_0(n)=\lim_n \frac{q_n(0)}{l_n^{\alpha}(0)}=b_0=\left\{
\begin{array}{ll}
0, & \hbox{if $M_0>0$;} \\
C\not=0, & \hbox{if $M_0=0$.}
 \end{array}
 \right.$$

\noindent Now, from (\ref{cocientederivadasortonormales}) for $k=1$, (\ref{limitecocientederivadas}) and (\ref{expresionA}) we get

$$ \lim_n b_1(n)=\lim_n \frac{1}{A_1(1,n)}\left(\frac{q'_n(0)}{(l_n^{\alpha})'(0)}-b_0(n)\right)=b_1$$
\noindent Observe that
$$b_1=\left\{
\begin{array}{ll}
0, & \hbox{if $M_0>0$ and $M_1>0$ ;} \\
C\not=0, & \hbox{if $M_0>0$ and $M_1=0$ .}
 \end{array}
 \right.$$
 \noindent In this way, recursively, if $M_0M_1 \dots M_i>0$ and $M_{i+1}=0$ we can assure that
\begin{equation*}
b_j=\left\{
\begin{array}{ll}
 0, & \hbox{if $0 \le j \le i$;} \\
 C\not=0, & \hbox{if $j=i+1$,}
 \end{array}\right.
\end{equation*}
and we obtain the result. $\quad \Box$

\bigskip

As a consequence of the  above lemma, we can establish a Mehler-Heine type formula for general discrete  Laguerre-Sobolev  orthonormal polynomials. This formula shows how the presence of the masses in  the discrete part of the inner product changes the asymptotic behavior around the origin. Moreover, it supplies information on the location and asymptotic distribution of the zeros of the polynomials in terms of the zeros of known special functions.

We recall the corresponding formula for orthonormal Laguerre polynomials (see \cite{sz})
\begin{equation}\label{MHLaguerre}
\lim_n\frac{l_n^{\alpha}(x/(n+k))}{n^{\alpha /2}}= \sqrt{\Gamma (\alpha +1)} x^{-\alpha /2} J_{\alpha}(2 \sqrt{x})
\end{equation}
uniformly on compact subsets of $\mathbb{C}$ and uniformly for $k \in \mathbb{N}\cup \{0 \}$, where $J_{\alpha}$ is the Bessel function of the first  kind.

\begin{proposition}\label{Mehler-Heine} The polynomials $q_n$ satisfy the following Mehler-Heine type formula:
\begin{equation}\label{MHLaguerre-Sobolev}
\lim_n\frac{q_n(x/n)}{n^{\alpha /2}}= \sqrt{\Gamma (\alpha +1)} \sum_{j=0}^{N+1} b_j\, x^{-\alpha /2} J_{\alpha +2j}(2 \sqrt{x})
\end{equation}
\noindent uniformly on compact subsets of $\mathbb{C}$.
\end{proposition}
\textbf{Proof.} The proof is a straightforward consequence of formula (\ref{formulaKoekoek}), Lemma \ref{limitecoeficientesKoekoek}  and (\ref{MHLaguerre}).   $\quad \Box$

\medskip

\noindent \textbf{Remark.} According to  Lemma \ref{limitecoeficientesKoekoek}, the first Bessel function which appears in (\ref{MHLaguerre-Sobolev}) corresponds with the first index $j$ such that $M_j=0$, in the inner product (\ref{innerproduct}). We want to highlight that
this  result  generalizes the one obtained in \cite[Theorem 3]{A-B-P-R}, where the authors only deal with  inner products with a unique $\lq\lq$gap" in the discrete part.

\bigskip
The above proposition allows us to deduce a lower estimate of $\Vert q_n \Vert _{L_p(x^{\beta}dx)}$, for $\beta=\alpha$ and $\beta=\alpha p/2$, that will play an important role in the proof of Theorem \ref{cohen}.

\medskip

\begin{proposition}\label{acotacionnormasSobolev} Let $1 \le p \le \infty$. Then, the following statements hold:

\noindent For $\alpha>-1/2$
\begin{equation*}
\Vert q_n \Vert _{L_p(x^{\alpha}dx)} \ge C \left\{
                                                  \begin{array}{ll}
                                                 n^{-1/4}(\log (n+1))^{1/p}, & \hbox{if $p=\frac{4\alpha+4}{2\alpha+1}$;} \\
                                                 n^{\alpha/2-(\alpha+1)/p}, & \hbox{if $\frac{4\alpha+4}{2\alpha+1}<p \le \infty ,$}
                                                  \end{array}
                                                  \right.
\end{equation*}
and for $\alpha>-2/p$ if $1 \le p < \infty$ and $\alpha \ge 0$ if $p = \infty$
\begin{equation*}
\Vert q_n \Vert _{L_p(x^{\alpha p/2}dx)} \ge C\left\{
\begin{array}{ll}
 n^{-1/4}(\log (n+1))^{1/p}, & \hbox{if $p=4$;} \\
 n^{-1/p}, & \hbox{if $4<p \le \infty $,}
 \end{array}
 \right.
\end{equation*}
where $C$ is an absolute positive constant.
\end{proposition}

\noindent \textbf{Proof.} Assume $1 \le p < \infty$. Then,
\begin{align*}
\nonumber
&\Vert q_n \Vert^p _{L_p(x^{\beta}dx)}= \int_0^{\infty}|q_n(x)e^{-x/2}|^px^{\beta}dx \\ \nonumber
&> \int_0^{1/\sqrt{n}}|q_n(x)e^{-x/2}|^px^{\beta}dx \ge C n^{- \beta-1}\int_0^{\sqrt{n}}|q_n(t/n)|^pt^{\beta}dt \\ \nonumber
\end{align*}
According to formula (\ref{MHLaguerre-Sobolev}), $\exists \, n_0\in \mathbb{N}$ such that $\forall n\ge n_0$
$$\int_0^{\sqrt{n}}|q_n(t/n)|^pt^{\beta}dt  \ge C n^{p\alpha/2} \int_0^{\sqrt{n}} | \sum_{j=0}^{N+1}  b_j\, t^{-\alpha/2}\,  J_{\alpha +2j}(2\sqrt{t})|^p t^{\beta} dt$$
and therefore $\forall n\ge n_0$
$$ \Vert q_n \Vert^p _{L_p(x^{\beta}dx)}
\ge Cn^{p\alpha/2- \beta-1}\int_0^{2{n}^{1/4}}u^{2\beta-p\alpha+1} | \sum_{j=0}^{N+1} b_j\,  J_{\alpha +2j}(u)|^pdu .$$
Working  as  Stempak in   \cite[Lemma 2.1]{S}, we can prove that
for $\alpha >-1$, and $\lambda>-1 -\alpha p$
\begin{equation*}
\int_0^{2{n}^{1/4}}u^{\lambda} | \sum_{j=0}^{N+1} b_j \,  J_{\alpha +2j}(u)|^pdu \sim \left\{
\begin{array}{ll}
1, & \hbox{if $\lambda<p/2-1$;} \\
\log (n+1), & \hbox{if $\lambda=p/2-1$.}
\end{array}
\right.
\end{equation*}
Thus, if $1 \le p < \infty$,  we obtain the first and the second result for $\beta=\alpha$ and  $\beta=p\alpha/2$ respectively. The results for $p=\infty$ can be deduced from the previous one by  passing to the limit when $p$ goes to $\infty$.
$\quad\Box$

\bigskip

It is worth to noticing that these lower bounds are sharp in the following sense.

\begin{proposition} Let $1 \le p \le \infty$. Then:

\noindent For $\alpha \ge 0$,
\begin{equation*}
\Vert q_n \Vert _{L_p(x^{\alpha}dx)} \sim \left\{
\begin{array}{ll}
n^{-1/4}(\log (n+1))^{1/p}, & \hbox{if $p=\frac{4\alpha+4}{2\alpha+1}$;} \\
n^{\alpha/2-(\alpha+1)/p}, & \hbox{if $\frac{4\alpha+4}{2\alpha+1}<p\le {\infty} $,}
\end{array}
\right.
\end{equation*}
and for  $\alpha>-2/p$ if $1 \le p < \infty$ and $\alpha \ge 0$ if $p = \infty$,
\begin{equation*}
\Vert q_n \Vert _{L_p(x^{\alpha p/2}dx)} \sim \left\{
\begin{array}{ll}
n^{-1/4}(\log (n+1))^{1/p}, & \hbox{if $p=4$;} \\
n^{-1/p}, & \hbox{if $4<p \le {\infty} $.}
\end{array}
\right.
 \end{equation*}
\end{proposition}

\noindent \textbf{Proof.} From Lemma 1 of \cite {M1} it can be deduced that
for $\alpha \ge 0$
\begin{equation*}
\int_0^{\infty}|x^jl_n^{\alpha +2j}(x)e^{-x/2}|^px^{\alpha}dx   \sim \left\{
\begin{array}{ll}
n^{-p/4}\log (n+1), & \hbox{if $p=\frac{4\alpha+4}{2\alpha+1}$;} \\
n^{\alpha p/2-(\alpha+1)}, & \hbox{if $\frac{4\alpha+4}{2\alpha+1}<p\le {\infty} $,}
\end{array}
\right.
\end{equation*}
and for  $\alpha>-2/p$ if $1 \le p < \infty$ and $\alpha \ge 0$ if $p = \infty$
\begin{equation*}
\int_0^{\infty}|x^jl_n^{\alpha +2j}(x)e^{-x/2}x^{\alpha/2}|^pdx  \sim \left\{
\begin{array}{ll}
n^{-p/4}\log (n+1), & \hbox{if $p=4$;} \\
n^{-1}, & \hbox{if $4<p \le {\infty} $.}
\end{array}
\right.
 \end{equation*}
 Thus, using  the representation formula  for  the polynomials $q_n $ (see (\ref{formulaKoekoek})),  and the fact that  the connection coefficients are bounded (see Lemma \ref{limitecoeficientesKoekoek}), we get one of the two inequalities. The other one has been proved in
Proposition \ref{acotacionnormasSobolev} and therefore  the result follows. $\quad\Box$

\section {A Cohen type inequality}

In this section we prove a Cohen type inequality for the Fourier expansions in terms of discrete Laguerre-Sobolev orthonormal polynomials with an arbitrary (finite) number of mass points. So
we extend the  result  due to Fejzullahu and  Marcell\'{a}n which deals with a discrete Laguerre-Sobolev inner product with at most  two masses  in the discrete part (see \cite{B-P}).
\bigskip

\noindent \textbf{Proof of Theorem \ref{cohen}.}  Let us consider the following test functions which were already used  in  \cite{M2} and later in \cite{B-P}
\begin{equation*}\label{funcionestest}
g_n^{\alpha,j}(x)=x^j\left[L_n^{\alpha+j}(x)-\sqrt{\frac{(n+1)(n+2)}{(n+\alpha+j+1)(n+\alpha+j+2)}}\,L_{n+2}^{\alpha+j}(x)
\right],
\end{equation*}
with  $j \in \mathbb{N}\setminus \{1, \dots, N \}$. Notice that

\begin{equation}\label{derivadasfuncionestest}
(g_n^{\alpha,j})^{(i)}(0)=0, \, i=0, \dots, N.
\end{equation}
These functions can be written as (see formula (2.15) in \cite{M2})

\begin{equation}\label{desarrollofuncionestest}
g_n^{\alpha,j}(x)=\sum_{m=0}^{j+2}a_{m,j}(\alpha, n)L_{n+m}^{\alpha}(x)
\end{equation}
with
\begin{equation*}\label{estimacionaoj}
a_{0,j}(\alpha, n)=\frac{\Gamma(n+\alpha+j+1)}{\Gamma(n+\alpha+1)}\cong n^j.
\end{equation*}
From (\ref{derivadasfuncionestest}), (\ref{desarrollofuncionestest}),  and $0 \le k \le  n$, we have
\begin{align*}
&\widehat{g_n^{\alpha,j}}(k)=\langle g_n^{\alpha,j}, q_k \rangle_S =\frac{1}{\Gamma(\alpha+1)}\int_0^{\infty}g_n^{\alpha,j}(x)q_k(x)e^{-x}x^{\alpha}dx \\
&=\frac{1}{\Gamma(\alpha+1)}\sum_{m=0}^{j+2}a_{m,j}(\alpha, n)\int_0^{\infty}L_{n+m}^{\alpha}(x)q_{k}(x)e^{-x}x^{\alpha}dx\,.
\end{align*}
By the orthogonality of Laguerre polynomials,  we obtain
\begin{equation*}
\widehat{g_n^{\alpha,j}}(k)=
\left\{
  \begin{array}{ll}
    0,  &\hbox{\, if $0 \le k \le n-1$;} \\
    \frac{1}{\Gamma(\alpha+1)}a_{0,j}(\alpha, n)\int_0^{\infty}L_{n}^{\alpha}(x)q_{n}(x)e^{-x}x^{\alpha}dx ,
     &\hbox{\, if $k=n$.}
  \end{array}
\right.
\end{equation*}
Thus,  from Lemma \ref{cocientederivadas} (b), the estimate of $a_{0,j}(\alpha, n)$ and the value of the norm of Laguerre polynomials (see (\ref{norma Laguerre})), we can deduce
\begin{align*}
& \widehat{g_n^{\alpha,j}}(n)=\frac{1}{\Gamma(\alpha+1)}a_{0,j}(\alpha, n)\int_0^{\infty}L_{n}^{\alpha}(x)\frac{Q_{n}(x)}{\langle Q_{n}, Q_{n} \rangle_S^{1/2} }e^{-x}x^{\alpha}dx =\\
&a_{0,j}(\alpha, n) \, \frac{\Vert L_n^{\alpha}\Vert^2}{\langle Q_{n}, Q_{n} \rangle_S^{1/2}} \cong a_{0,j}(\alpha, n) \, \Vert L_n^{\alpha}\Vert  \cong  \frac{n^{j +\alpha/2}}{\sqrt{\Gamma(\alpha+1)}}\\
\end{align*}
Observe that $Q_{n}$ and $L_{n}^{\alpha}$ have always equivalent norms, and, therefore this estimation does not depend neither on the number of positive masses, nor on the existence or non-existence of any gap in the inner product.

Applying the operator $T_n^{\alpha,S}$ to the functions $g_n^{\alpha,j}$, we get
\begin{equation*}\label{Tfunciones test}
T_n^{\alpha,S}(g_n^{\alpha,j})= c_{n,n}\, \widehat{g_n^{\alpha,j}}(n)q_n,
\end{equation*}
and therefore
\begin{align*}
\Vert T_n^{\alpha,S}\Vert_{[S_p^{\beta}]} & \ge (\Vert g_n^{\alpha,j}\Vert_{S_p^{\beta}})^{-1}\Vert T_n^{\alpha,S}(g_n^{\alpha,j}) \Vert_{S_p^{\beta}}=(\Vert g_n^{\alpha,j}\Vert_{S_p^{\beta}})^{-1} |c_{n,n}||\widehat{g_n^{\alpha,j}}(n)|\,\Vert q_n \Vert_{S_p^{\beta}}\\
&\ge (\Vert g_n^{\alpha,j}\Vert_{S_p^{\beta}})^{-1} |c_{n,n}||\widehat{g_n^{\alpha,j}}(n)|\,\Vert q_n \Vert_{L_p(x^{\beta}dx)}.\\
\end{align*}
\noindent On the other hand, for $j>\alpha-1/2-2(\alpha +1)/p$ we have

\begin{equation*}
\Vert g_n^{\alpha,j}\Vert_{S_p^{\beta}} \le c \left\{
                                                  \begin{array}{ll}
                                                    n^{j-1/2+(\alpha +1)/p}, & \hbox{if $\beta=\alpha$;} \\
                                                    n^{\alpha/2+j-1/2+1/p}, & \hbox{if $\beta=p\alpha/2$,}
                                                  \end{array}
                                                \right.
\end{equation*}
(see formula (3.3) and formula (1.19), (2.12) in \cite{M2} respectively).
 Thus,  by Proposition \ref{acotacionnormasSobolev} we get:

\noindent For $\beta=\alpha$ with $\alpha >-1/2$
\begin{equation*}
\Vert T_n^{\alpha,S}\Vert_{[S_p^{\alpha}]} \ge C |c_{n,n}|\left\{
                                                   \begin{array}{ll}
                                                     (\log (n+1))^\frac{2\alpha +1}{4\alpha +4}, & \hbox{if $p=q_0$;} \\
                                                     n^{\alpha+1/2-2(\alpha+1)/p}, & \hbox{if $q_0<p \le \infty$.}
                                                   \end{array}
 \right.
\end{equation*}

\noindent For $\beta=p\alpha/2$ with $\alpha>-2/p$ if $1 \le p<\infty$ and $\alpha\ge 0$ if $p=\infty$,
\begin{equation*}
\Vert T_n^{\alpha,S}\Vert_{[S_p^{p\alpha/2}]} \ge C|c_{n,n}|\left\{
                                                   \begin{array}{ll}
                                                     (\log (n+1))^{1/4}, & \hbox{if $p=4$;} \\
                                                     n^{1/2-2/p}, & \hbox{if $4<p \le \infty$.}
                                                   \end{array}
\right.
\end{equation*}
 Hence, by duality the theorem follows. $\quad \Box$

\medskip

\noindent\textbf{Remark.} In particular, for $M_i=0$, $i=0, \dots, N$, the above theorem extends  Theorem 1 in \cite{M2} to negative values of $\alpha$.

\bigskip

In the particular case of $c_{k,n}=1$, $k=0, \dots, n$, the operator $T_n^{\alpha,S}$ is the  $n$th partial sum of the Fourier expansion, so, we can assure the   following result.

\begin{corollary} If $p$ is outside the Pollard interval $(p_0, q_0)$, we have
\begin{equation*}
\Vert S_n \Vert_{[S_p^{\beta}]} \rightarrow \infty,\quad  n \rightarrow \infty
\end{equation*}
where $S_n$ denotes the $n$th partial sum of the Fourier expansion.
\end{corollary}

\end{document}